%% file: Joint-Surv-Sign-rev-200811.tex
\documentclass[preprint,12pt]{elsarticle}

\makeatletter
\def\ps@pprintTitle{%
 \let\@oddhead\@empty
 \let\@evenhead\@empty
 \def\@oddfoot{\centerline{\thepage}}%
 \let\@evenfoot\@oddfoot}
\makeatother

\usepackage{anysize} 
\marginsize{1in}{1in}{.9in}{1in}
\usepackage{hyperref}
\usepackage{multirow}

\usepackage{epsfig}
\usepackage{amsmath,amssymb, amstext,mathtext,amsthm,color}

\usepackage{amssymb}

\usepackage{graphicx,color}
\usepackage[numbers]{natbib}

\usepackage{amsmath,amssymb}

\usepackage{comment}
\usepackage{float}
\usepackage{setspace} 
\usepackage{epstopdf}

\newtheorem{example}{Example}

\usepackage{multirow,xcolor,colortbl,rotating}

\usepackage{xcolor}     
\usepackage{mdframed}   

\begin{document}
\begin{frontmatter}
\title{The joint survival signature  of  coherent systems with shared components}
\author{Tahani Coolen-Maturi$^{1,}$\footnote{Corresponding author: tahani.maturi@durham.ac.uk}, Frank P.A.\ Coolen$^{1}$, 
Narayanaswamy  Balakrishnan$^{2}$}
\address{$^{1}$Department of Mathematical Sciences, Durham University, Durham, UK.\\ 
$^{2}$Department of Mathematics and Statistics,
McMaster University, Ontario, Canada.}

\begin{abstract} 
The concept of joint bivariate signature, introduced by \citet{navarro2013}, is a useful tool for quantifying the
reliability of two systems with shared components. As with the univariate system signature, introduced by Samaniego 
\cite{SAM07}, its applications are limited to systems with only one type of components, which restricts its  practical 
use. Coolen and Coolen-Maturi \cite{CCM12} introduced the survival signature, which generalizes Samaniego's signature
and can be used for systems with multiple types of components. This paper introduces a joint survival signature for multiple 
systems with multiple types of components and with some components shared between systems. A particularly important
feature is that the functioning of these systems can be considered at different times, enabling computation of relevant
conditional probabilities with regard to a system's functioning conditional on the status of another system with which it
shares components. Several opportunities for practical application and related challenges for further development of the 
presented concept are briefly discussed, setting out an important direction for future research.

\end{abstract}

\begin{keyword}
Coherent systems \sep exchangeable components \sep  signature \sep survival signature \sep system reliability.
\end{keyword}
\end{frontmatter}

\section{Introduction}
\label{intro}
In recent decades, the system signature has become a popular tool for quantifying reliability of coherent systems consisting 
of components with exchangeable random failure times \cite{SAM85}, where in the literature the assumption of exchangeability
\cite{DF74} is often replaced by the stronger assumption of independent and identically distributed ({\em iid}) component 
failure times. The system signature is a summary of the system structure function which is sufficient to quantify several
important aspects of reliability of a system, in particular the system's failure time distribution. A detailed introduction 
and overview of system signatures is presented by Samaniego \cite{SAM07}.

The essential property of the system signature is that it enables information of the system structure to be fully taken 
into account through the signature, and this is separated from information about the random failure times of the components. 
The main disadvantage of system signatures, however, is that it becomes extremely complicated, and is indeed effectively 
impossible, to keep this separation when generalizing the concept to systems with multiple types of components, which is 
crucial for a practically applicable theory as most real-world systems consist of more than a single type of components 
\cite{CCM12,NSB11}. As an alternative to the system signature, Coolen and Coolen-Maturi \cite{CCM12} introduced the 
survival signature. For systems with just one type of components, the survival signature is equivalent to 
the system signature, but the survival signature can be defined for, and easily applied to, systems with multiple 
types of components.  

There are many scenarios where two or more systems share components, which can be of different types. Consider for 
example the case of two computers linked to a server, where the performance of any computer will depend on the 
performance of the shared components, in the server, and the performance of its own components. It should be emphasized
that systems has a wide meaning in this context, including not only engineering systems but also networks and organisational
structures. For example, if good functioning of multiple academic departments at one university during an exams period
with strict marking deadlines depend on one central information technology support group, then the latter can be regarded
as a component shared by the different departments (`systems'). While we do not focus explicitly on it, it is important 
to note that the theory in this paper can also be applied for the case of one system which performs two or more functions, 
with some but not all components involved in multiple functions.

\citet{NSB10} and \citet{ZMB17} introduced the system signature for systems with shared components, however the joint 
system signature representation has no direct probability interpretation as it can take negative values. 
\citet{navarro2013} presented the so-called joint bivariate system signature, which has a probabilistic interpretation, 
and they also showed how the joint bivariate system signature can be used to perform stochastic comparisons. But 
again their method is limited to one type of component, which is less useful in real world applications.  

In this paper, we introduce the joint survival signature of coherent systems with shared components, which can be 
of different types. This first presentation of the new concept emphasizes the opportunity to consider the reliability 
of the systems at different time points, which is crucial for many practical applications. In particular, it enables one
to infer one system's reliability conditional on the information that another system, with which it shares some components,
functions or not at a different time point. This introduction of the new concept is the first step of an important and
extensive research direction, challenges for computation and application will be outlined in the final section.

This paper is organised as follows. Section \ref{sec.sur.sign} provides a brief overview of the survival signature 
introduced by Coolen and Coolen-Maturi \cite{CCM12}. Section \ref{sec.twosign} introduces the new joint survival 
signature of two coherent systems with shared components, followed by generalisation to three coherent 
systems with shared components in Section \ref{sec.three}. It is briefly discussed how the further generalization 
to more than three systems can be achieved, all required ingredients for such a generalization follow quite straightforwardly
from the case with three systems. Finally, in Section \ref{concl} opportunities and challenges for further research to 
enable practical application of the new concept to large-scale systems and networks are briefly discussed.

\section{Survival signature}
\label{sec.sur.sign}

For a system with $n$ components, we define the state vector $\underline{x} \in \{0,1\}^n$ with entry
$x_i=1$ if component $i$ functions and $x_i=0$ if not. The labelling of the components is arbitrary but must be 
fixed to define $\underline{x}$. The structure function $\phi : \{0,1\}^n \rightarrow \{0,1\}$, defined for all possible 
$\underline{x}$, takes the value 1 if the system functions and 0 if the system does not function for state vector 
$\underline{x}$. In this paper, we restrict attention to coherent systems, which implies that $\phi(\underline{x})$ 
is not decreasing in any of the components of $\underline{x}$, so system functioning cannot be improved by worse 
performance of one or more of its components. We further assume that $\phi(\underline{0})=0$ and $\phi(\underline{1})=1$, 
so the system fails if all its components fail and it functions if all its components function. These assumptions could be 
relaxed but are reasonable for most practical systems.

Consider a system with $K\geq 2$ types of components, with $n_k$ components of type $k \in \{1,2,\ldots,K\}$ and 
$\sum_{k=1}^K n_k = n$. Assume that the random failure times of components of the same type are exchangeable \cite{DF74}, 
while full independence is assumed for the random failure times of components of different types. 
Due to the arbitrary ordering of the components in the state vector, 
components of the same type can be grouped together, leading to a state vector that can be written as 
$\underline{x} = (\underline{x}^1,\underline{x}^2,\ldots,\underline{x}^K)$, with 
$\underline{x}^k = (x^k_1,x^k_2,\ldots,x^k_{n_k})$ the sub-vector representing the states of the components of 
type $k$. 

Coolen and Coolen-Maturi \cite{CCM12} introduced the {\em survival signature} for such a system, denoted by 
$\Phi(l_1,l_2,\ldots,l_K)$, with $l_k=0,1,\ldots,n_k$ for $k=1,\ldots,K$, which is defined to be the probability that the system 
functions given that {\it precisely} $l_k$ of its $n_k$ components of type $k$ function, for each $k\in \{1,2,\ldots,K\}$. 

There are $\binom{n_k}{l_k}$ state vectors $\underline{x}^k$ with $\sum_{i=1}^{n_k} x^k_i = l_k$; let $S^k_l$ 
denote the set of these state vectors for components of type $k$ and let $S_{l_1,\ldots,l_K}$ denote the set of 
all state vectors for the whole system for which $\sum_{i=1}^{n_k} x^k_i = l_k$, $k=1,2,\ldots,K$. Due to the 
exchangeability assumption for the failure times of the $n_k$ components of type $k$, all the state vectors 
$\underline{x}^k \in S^k_l$ are equally likely to occur, hence
\begin{equation} \label{phi-mult}
\Phi(l_1,\ldots,l_K) = \left[ \prod_{k=1}^K \binom{n_k}{l_k}^{-1} \right] \times 
\sum_{\underline{x} \in S_{l_1,\ldots,l_K}} \phi(\underline{x})
\end{equation}

Let $C_k(t) \in \{0,1,\ldots,n_k\}$ denote the number of components of type $k$ in the system which function 
at time $t>0$. The probability that the system functions at time $t>0$ is 
\begin{equation} 
P(T_S>t) = \sum_{l_1=0}^{n_1} \cdots \sum_{l_K=0}^{n_K}  \Phi(l_1,\ldots,l_K) P(\bigcap_{k=1}^K \{C_k(t) = l_k\})  \label{surv-sys-mult}
\end{equation}
If one assumes independence of the failure times of components of different types, then this leads to, for $l_k \in \{0,1,\ldots,n_k\}$ for each $k\in \{1,\ldots,K\}$,
\begin{equation} 
P(\bigcap_{k=1}^K \{C_k(t) = l_k\}) = \prod_{k=1}^K P(C_k(t) = l_k) 
\end{equation}
If, in addition, one assumes that the failure times of components of the same type are {\em iid} with known 
cumulative distribution function (CDF) $F_k(t)$ for type $k$, then this leads to
\begin{equation} 
P(\bigcap_{k=1}^K \{C_k(t) = l_k\}) = \prod_{k=1}^K \binom{n_k}{l_k} [F_k(t)]^{n_k-l_k} [1-F_k(t)]^{l_k} 
\end{equation}

A crucial practical consideration is how to decide if it is reasonable to assume that components are of the same type, in the
sense of having exchangeable failure times. An easy way to think about this is as follows. Suppose that there is a number of components,
and you get the information that, at a specific time, one of them has failed, without any further information which enables
you to identify which component has failed. Exchangeability then implies that each of these components is equally likely to be
the one that has failed. Note that this includes consideration of the role of the components and the environment in which they
function in the system. Of course, this is a subjective modelling assumption which relates to the level of detail in which one
models the system, in most of the reliability theory literature it is silently assumed. In particular the common assumption of
{\em iid} failure times of components is a stronger assumption.

Since its presentation by Coolen and Coolen-Maturi \cite{CCM12}, there has been substantial research contributing to the
further theory and applicability of survival signature methods. An important topic is computation of the survival signature,
a useful method based on binary decision diagrams has been presented by Reed \cite{Re17}, while derivation of the survival
signature for systems built up by subsystems in series or parallel configuration was also presented \cite{Co14}. The survival
signature also enables very efficient simulation methods to be developed \cite{Ge19,PATELLI2017327}, and further examples of powerful
methodology for system reliability quantification enabled by the use of survival signatures include the modelling of dependence 
between components of different types \cite{ERYILMAZ2018118,Ge19}, Bayesian and nonparametric predictive inference \cite{As15,Co14}, 
reliability-redundancy allocation \cite{HUANG2019511}, phased-missions \cite{Hu19}, component reliability importance 
measures \cite{Fe16}, resilience achieved by swapping components within a system \cite{Na19} and stochastic comparison of 
different systems \cite{Sa16}.

\section{Joint survival signature of two coherent systems with shared components}
\label{sec.twosign}

In this section we present the joint survival signature of two coherent systems that share some components. It is
important that the functioning of the two systems can be considered at different moments in time. This can be used 
to derive the marginal survival function of one of the systems as well as the conditional reliability of a system given
the status of the other system at any time. First we consider systems that only have a single type of components, which
is helpful to explain the main ideas and notation. This is later extended to systems with components of multiple types.
It should be emphasized that no maintenance or replacement activities are being considered throughout this paper, so once
a component has failed it remains in failed state.

\subsection{One-type of components} 
\label{subsec.onetype}

Let $T_1$ and $T_2$ be the failure times of two coherent systems, $S_1$ and $S_2$, based on components with {\em iid} failure times $X_1, \ldots,X_n$ having a common continuous distribution function $F$. The assumption of {\em iid} failure times is for ease of presentation. 
 A coherent system fails at the failure of one of its components. Assume that the first system has $n^*_1$ components and the second system has 
$n^*_2$ components, with $n_{12}$ of these components in common, these are called shared components, so in total there are $n=n^*_1+n^*_2-n_{12}$ components in the two systems. The case of interest in this paper is $n_{12}>0$, otherwise the two systems are independent by the {\it iid} assumption. Let the numbers of components in $S_1$ and $S_2$ which are not shared with the other system be denoted by $n_1$ and
$n_2$, respectively, so $n^*_1=n_1+n_{12}$, $n^*_2=n_2+n_{12}$ and $n=n_1+n_2+n_{12}$.

The joint survival signature will enable reliability quantification for both systems at possibly different times, say $S_1$ is considered
at time $t_1$ and $S_2$ at time $t_2$. This means that the numbers of the shared components functioning at these two different times
must both be specified, note that the specific times $t_1$ and $t_2$ do not play any further role in the survival signature, as this
enables inference for all time points when combined with the component failure time distributions. This means that the joint survival
signature presented in this paper has the same advantageous property as the survival signature for a single system, that is it takes
the structure of the systems into account while being separated from the random component failure times.

The joint survival signature $\Phi(l_1,l_2,l_{[1]2},l_{1[2]})$ can be defined as the probability that systems $S_1$ and $S_2$ both function 
given that precisely $l_1$ out of $n_1$ and $l_{[1]2}$ out of the $n_{12}$ shared components function when $S_1$ is being considered,
and precisely $l_2$ out of $n_2$ and $l_{1[2]}$ out of the $n_{12}$ shared components function when $S_2$ is being considered. Denoting
the events that $S_1$ and $S_2$ function at the moment of time they are considered by $SF_1$ and $SF_2$, respectively, the joint
survival signature denoted by

\begin{equation}
\Phi(l_1,l_2,l_{[1]2},l_{1[2]})= P(SF_1,  SF_2  | \\
 l_1 , l_2 , l_{[1]2}, l_{1[2]}  )
\end{equation}

It is also important to emphasise, given the above setting,  that the same $\min(l_{1[2]},l_{[1]2})$ shared components are functioning at both times $t_1$ and $t_2$, and the same $n_{12} - \max(l_{1[2]},l_{[1]2})$ shared components are not functioning at both times. The remaining components (if $l_{1[2]}\neq l_{[1]2}$ are different) fail between the two different times. Therefore, $\Phi(l_1,l_2,l_{[1]2},l_{1[2]}) =0$  if  $t_1<t_2$ and $l_{1[2]} \geq l_{[1]2}$, or if $t_1>t_2$ and $l_{1[2]}\leq l_{[1]2}$.

Let $C^1_{t_1}\in \{0,1,\ldots,n_1\}$ denote the number of components that are only in system $S_1$ that function at time $t_1>0$, and
let $C^{[1]2}_{t_1}\in \{0,1,\ldots,n_{12}\}$ denote the number of shared components in system $S_1$ that function at $t_1$. 
Similarly, let $C^2_{t_2}\in \{0,1,\ldots,n_2\}$ denote the number of components that are only in system $S_2$ that function at time 
$t_2>0$, and let $C^{1[2]}_{t_2}\in \{0,1,\ldots,n_{12}\}$ denote the number of shared components in system $S_2$  that function at $t_2$. 
Let the probability distribution of the component failure time have CDF $F(t)$, then 
for $t_1<t_2$, which implies that $l_{1[2]} \leq l_{[1]2}$, 
\begin{align}
P(T_1>t_1,T_{2}>t_2)=\sum_{l_1=0}^{n_1}\sum_{l_2=0}^{n_2}\sum_{l_{[1]2}=0}^{n_{12}}\sum_{l_{1[2]}=0}^{n_{12}}\Phi(l_1,l_2,l_{[1]2},l_{1[2]}) P_{C_{t_1<t_2}}(l_1,l_2,l_{[1]2},l_{1[2]})
\end{align}
where the assumption of {\em iid} component failure times leads to
\begin{align}
\nonumber &P_{C_{t_1<t_2}}(l_1,l_2,l_{[1]2},l_{1[2]})=P(C^1_{t_1}=l_1,C^2_{t_2}=l_2,C^{[1]2}_{t_1}=l_{[1]2},C^{1[2]}_{t_2}=l_{1[2]})\\
\nonumber &=\frac{n_{1}!}{(n_{1}-l_{1})! l_1!} [1-F(t_1)]^{l_1}[F(t_1)]^{n_1-l_1} \frac{n_{2}!}{(n_{2}-l_{2})! l_2!} [1-F(t_2)]^{l_2}[F(t_2)]^{n_2-l_2}\\
&\times \frac{n_{12}!}{(n_{12}-l_{[1]2})!  (l_{[1]2}-l_{1[2]})! l_{1[2]}!}[F(t_1)]^{n_{12}-l_{[1]2}}[F(t_2)-F(t_1)]^{l_{[1]2}-l_{1[2]}}[1-F(t_2)]^{l_{1[2]}}
\end{align}

\noindent Similarly, for $t_1>t_2$ which implies that $l_{1[2]} \geq l_{[1]2}$, 
\begin{align}
P(T_1>t_1,T_{2}>t_2)=\sum_{l_1=0}^{n_1}\sum_{l_2=0}^{n_2}\sum_{l_{[1]2}=0}^{n_{12}}\sum_{l_{1[2]}=0}^{n_{12}}\Phi(l_1,l_2,l_{[1]2},l_{1[2]}) P_{C_{t_1>t_2}}(l_1,l_2,l_{[1]2},l_{1[2]})
\end{align}
where, under the {\em iid} assumption,
\begin{align}
\nonumber &P_{C_{t_1>t_2}}(l_1,l_2,l_{[1]2},l_{1[2]})=\\
\nonumber &\frac{n_{1}!}{(n_{1}-l_{1})! l_1!} [1-F(t_1)]^{l_1}[F(t_1)]^{n_1-l_1}\frac{n_{2}!}{(n_{2}-l_{2})! l_2!}[1-F(t_2)]^{l_2}[F(t_2)]^{n_2-l_2}\\
&\times \frac{n_{12}!}{(n_{12}-l_{1[2]})!  (l_{1[2]}-l_{[1]2})! l_{[1]2}!}[F(t_2)]^{n_{12}-l_{1[2]}}[F(t_1)-F(t_2)]^{l_{1[2]}-l_{[1]2}}[1-F(t_1)]^{l_{[1]2}}
\end{align}

If $t_1=t_2=t$, then, with notation $l_{12}=l_{[1]2}=l_{1[2]}$ and the {\em iid} assumption, 
\begin{align}
\nonumber &P_{C_{t_1=t_2}}(l_1,l_2,l_{[1]2},l_{1[2]})=\\
\nonumber &\frac{n_{1}!}{(n_{1}-l_{1})! l_1!} [1-F(t_1)]^{l_1}[F(t_1)]^{n_1-l_1} \frac{n_{2}!}{(n_{2}-l_{2})! l_2!} [1-F(t_2)]^{l_2}[F(t_2)]^{n_2-l_2}\\
&\times \frac{n_{12}!}{(n_{12}-l_{12})!   l_{12}!}[F(t_1)]^{n_{12}-l_{12}}[1-F(t_2)]^{l_{12}}
\end{align}

These joint survival functions for systems $S_1$ and $S_2$ provide a detailed quantification of the reliability
of the two systems, based on the probability distributions for the component failure times. Note that, indeed, the structural
aspects of the systems are all taken into account by the newly proposed joint survival signature, which is independent of time,
while all temporal aspects are taken into account by the failure time distributions, and these two important aspects are fully
separated. This joint survival signature provides exciting opportunities for study of a variety of theoretical aspects and 
applications, in line with the contributions to the literature for the single system survival signature mentioned above. 

The following two examples illustrate the above introduced joint survival signature. Example \ref{Ex1} shows in detail how
the value of the joint survival signature can be derived for a single specific input. Example \ref{Ex2} presents the joint
survival functions for two basic systems with shared components.

\begin{example}
\label{Ex1}{\rm

Consider the two systems in Figure \ref{fig.ex1a}. All components are assumed to be of the same type, with exchangeable
failure times. The systems share components A, B and C, and each system has two further components. For System 1 we have 
$n^*_1=5$, $n_1=2$, $n_{12}=3$ and for System 2 we have $n^*_2=5$, $n_2=2$, $n_{12}=3$. For input ($l_1=1, l_2=1, l_{[1]2}=2, l_{1[2]}=1$),  Table \ref{tab.ex1a} lists all 24 possible scenarios of functioning components. If these are indeed the numbers of functioning
components, then each of these scenarios has probability 1/24 due to the exchangeability assumption. In this case, the system
functions for 10 of the 24 possibille scenarios, hence the survival signature is equal to $\Phi(1,1,2,1)=\frac{10}{24}$.

\input{fig.ex1a}

\begin{table}[H]
  \centering\small
    \begin{tabular}{ccccc}
    \multicolumn{1}{l}{System functions} & \multicolumn{1}{r}{$l_1=1$} & \multicolumn{1}{r}{$l_2=1$}  & \multicolumn{1}{r}{$l_{[1]2}=2$} & \multicolumn{1}{r}{$l_{1[2]}=1$} \\
\hline
    0     & D    & F  & AB         & A \\
    0     & D      & F& AB         & B \\
    0     & D    & G  & AB        & A \\
    0     & D   & G  & AB         & B \\
    0     & D   & F  & AC         & A \\
    1     & D   & F  & AC         & C \\
    1     & D   & G  & AC         & A \\
    0     & D    & G  & AC        & C \\
    1     & D     & F & BC        & B \\
    1     & D     & F& BC         & C \\
    0     & D       & G& BC       & B \\
    0     & D    & G  & BC        & C \\
    0     & E     & F& AB         & A \\
    1     & E     & F& AB         & B \\
    1     & E     & G & AB        & A \\
    0     & E   & G  & AB         & B \\
    0     & E    & F  & AC        & A \\
    1     & E    & F  & AC        & C \\
    1     & E    & G  & AC        & A \\
    0     & E     & G& AC         & C \\
    1     & E     & F& BC         & B \\
    1     & E     & F& BC         & C \\
    0     & E   & G   & BC        & B \\
    0     & E    & G  & BC        & C \\
\hline
    \end{tabular}
  \label{tab.ex1a}
\caption{Functioning components and system state, Example \ref{Ex1}}
\end{table}

}\end{example}

\begin{figure}[htbp]
\begin{center}
\includegraphics[scale=.3]{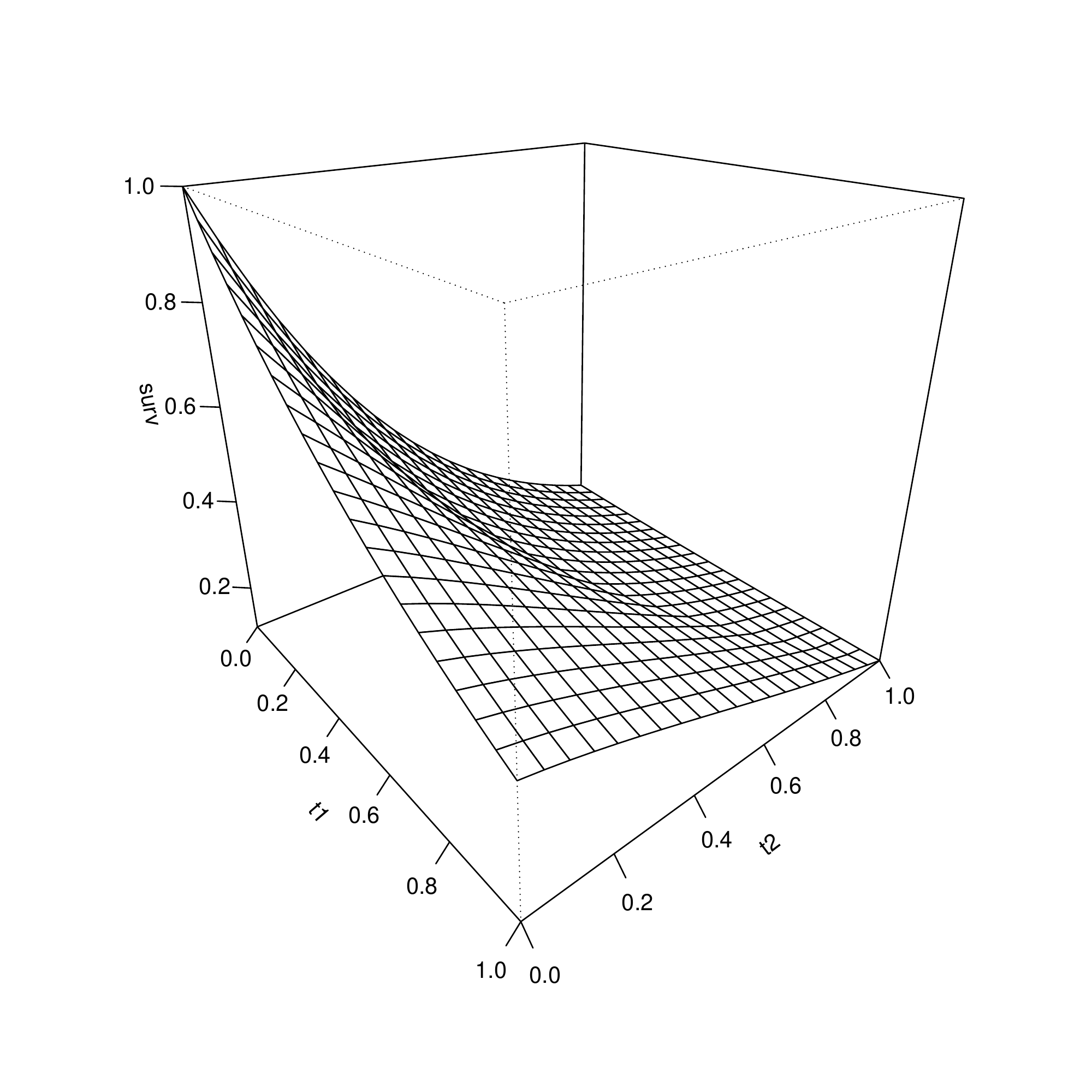}
\caption{Joint survival function of two systems  with one type of components, Example \ref{Ex2}}
\label{fig2}
\end{center}
\end{figure}

\begin{example}\label{Ex2}{\rm
Consider the two systems in Figure \ref{fig.ex2}, which share components A and B and each have one further components.
Again, all components are of the same type so their failure times are assumed to be exchangeable.  For System 1 we 
have $n^*_1=3$, $n_1=1$, $n_{12}=2$ and for System 2 we have $n^*_2=3$, $n_2=1$, $n_{12}=2$. The survival signature is 
zero, $\Phi(l_1,l_2,l_{[1]2},l_{1[2]})=0$, for the trivial cases: $l_1+l_2<2$, $l_{[1]2}=0$ and $l_{1[2]}=0$. The survival signature 
for remaining cases is given in Table \ref{tab.ex2}, it is derived in a similar way as illustrated in Example \ref{Ex1} for each possible
input, although the actual work required is of course not as bad as it may seem due to logical relationships between the systems'
states and different inputs under the assumption that the systems are coherent. Assuming that the components' failure times are {\em iid}
with Exponential distribution with rate 1, that is  $F(t)=1-e^{- t}$, then the joint survival function for the two system failure times 
is presented in Figure \ref{fig2}.

\input{fig.ex2} 
\begin{table}[H]
  \centering\small
    \begin{tabular}{ccccc}
       \multicolumn{1}{c}{C} & \multicolumn{1}{c}{D} & \multicolumn{2}{c}{A,B} &  \\
    \multicolumn{1}{c}{$l_1\in \{0,1\}$} & \multicolumn{1}{c}{$l_2\in \{0,1\}$} & \multicolumn{1}{c}{$l_{ [1] 2}\in \{0,1,2\}$} 
		&  \multicolumn{1}{c}{$l_{1[2]}\in \{0,1,2\}$} & \multicolumn{1}{c}{$\Phi(l_1,l_2,l_{[1]2},l_{1[2]})$} \\[1ex]
\hline
0	&	1	&	1	&	1	&	$1/2$	\\
1	&	0	&	1	&	1	&	0	\\
1	&	1	&	1	&	1	&	$1/2$	\\
0	&	0	&	1	&	1	&	0	\\
0	&	0	&	2	&	1	&	0	\\
0	&	1	&	2	&	1	&	$1/2$	\\
1	&	0	&	2	&	1	&	0	\\
1	&	1	&	2	&	1	&	$1/2$	\\
0	&	0	&	1	&	2	&	$1/2$	\\
0	&	1	&	1	&	2	&	$1/2$	\\
1	&	0	&	1	&	2	&	1	\\
1	&	1	&	1	&	2	&	1	\\
1	&	1	&	2	&	2	&	1	\\
\hline
    \end{tabular}
  \label{tab.ex2}
\caption{Two systems with one type of components, Example \ref{Ex2}}
\end{table}
}\end{example}

\subsection{Marginal and conditional survival functions}
\label{subsec.marginal}

If one has the joint survival signature for the two systems available, it is straightforward to derive the
marginal survival distribution for one of the systems, using the assumption that all components function
at time $0$, so
\begin{equation}
P(T_2>t_2)=P(T_1>0,T_{2}>t_2)=\sum_{l_2=0}^{n_2}\sum_{l_{1[2]}=0}^{n_{12}}\Phi(n_1,l_2,n_{12},l_{1[2]}) 
P_{C_{t_1<t_2}}(n_1,l_2,n_{12},l_{1[2]})
\end{equation}
where
\begin{equation}
P_{C_{t_1<t_2}}(n_1,l_2,n_{12},l_{1[2]})=
 \binom{n_2}{l_2}[1-F(t_2)]^{l_2}[F(t_2)]^{n_2-l_2} \binom{n_{12}}{l_{1[2]}}[F(t_2)]^{n_{12}-l_{1[2]}}[1-F(t_2)]^{l_{1[2]}}
\end{equation}

The joint survival signature can also be used to derive the conditional survival function of one system given the status, at any moment
of time, of the other system, with which it shares some components. Suppose that system $S_2$ is known to function at time $t_2$, but no 
further information is known about specific components. Conditional on this information, the probability that system $S_1$ functions 
at time $t_1$ is
\begin{align}
P(T_1>t_1|T_2>t_2)=\frac{P(T_1>t_1,T_{2}>t_2)}{P(T_2>t_2)}
\end{align}
Note that there are no limitations on the values of $t_1$ and $t_2$. While from theoretical perspective it is good to have the flexibility
to investigate the conditional survival functions for all values of $t_1$ and $t_2$, in practice it seems most natural that this probability
will be relevant for $t_1\geq t_2$.
  
It is also important to derive the conditional probability that system $S_1$ still
functions at time $t_1$ given the information that system $S_2$ has failed by time $t_2$. This is equal to
\begin{align}
P(T_1>t_1|T_2\leq t_2)=\frac{P(T_1>t_1,T_{2}\leq t_2)}{P(T_2\leq t_2)}
\end{align}
The probability in the numerator on the right-hand side of this equation has not been presented explicitly before in this paper, but
it can be derived in a similar way as the joint survival probability, presented above, by splitting into a new joint (survival) signature and
a factor related to the component failure time distributions. This requires a different (survival) signature, namely the probability for the
event that $S_1$ functions when considered, and $S_2$ does not function when considered; this is left as an exercise for the reader,
the steps involved follow the presentation above closely.

Also joint probabilities of system functioning at specific times, conditional on information of one or both of the systems functioning 
at other times, can be derived similarly. For example, if system 2 is known to function at $t_2$, the conditional probability of the
event that both systems function at $t>t_2$ is
\begin{align}
P(T_1>t,T_2>t|T_2>t_2)=P(\min\{T_1,T_2\}>t|T_2>t_2)=\frac{P(T_1>t,T_{2}>t)}{P(T_2>t_2)}
\end{align}
All such conditional probabilities are easily derived using the results presented 
above, or quite straightforward extensions of such results if variations to the joint survival signature are required.

\subsection{Multiple-types of shared components} 
\label{subsec.multiple}

The major advantage of the survival signature over Samaniego's system signature is the ability to use it for systems with
multiple types of components. The presentation above was restricted to a single type of component, in order to clearly introduce
the new concept of joint survival signature in detail. Now we indicate how the concept and results of the previous sections 
can be generalized to multiple types of components. This generalization is conceptually straightforward so it is only
briefly presented. We assume that there are $K$ component types in two systems, and denote the type of component 
$k \in \{1,2,\ldots,K\}$ as superscript. Let $C^{k,1}_{t_1}\in \{0,1,\ldots,n^k_1\}$, 
$C^{k,2}_{t_2}\in \{0,1,\ldots,n^k_2\}$, $C^{k,[1]2}_{t_1}\in \{0,1,\ldots,n^k_{12}\}$  
and $C^{k,1[2]}_{t_2}\in \{0,1,\ldots,n^k_{12}\}$ denote the numbers of components of type $k$, $k \in \{1,2,\ldots,K\}$,  
in system $S_1$ and system $S_2$  that function at time $t_1>0$ and $t_2>0$, respectively. For ease of notation, 
let $\underline{l}^k=(l^1_1,l^1_2,l^1_{[1]2},l^1_{1[2]},\ldots,l^K_1,l^K_2,l^K_{[1]2},l^K_{1[2]})$, and we denote the 
summation over all their possible values by $\sum_{\underline{l}^k}$. The joint survival signature for this scenario is
the probability that both systems function at the time points at which they are considered, given that the numbers of 
functioning components are represented by $\underline{l}^k$. If the probability distribution of the 
failure times of components of type $k$ has CDF $F_k(.)$, and assuming independence of failure times of components of
different types, we have 
\begin{equation}
P(T_1>t_1,T_{2}>t_2)=\sum_{\underline{l}^k}\Phi(\underline{l}^k) P(C_{\underline{t},\underline{l}^k})
\end{equation}
The detailed form of probability $P(C_{\underline{t},\underline{l}^k})$ depends again on the values of $t_1$ and $t_2$, 
for example for $t_1<t_2$ we have

\begin{align}
\nonumber P(C_{\underline{t},\underline{l}^k})&=\prod_{k=1}^K P(C^{k,1}_{t_1}=l^k_1,C^{k,2}_{t_2}=l^k_2,C^{k,[1]2}_{t_1}=l^k_{[1]2},C^{k,1[2]}_{t_2}=l^k_{1[2]})\\
\nonumber &=\prod_{k=1}^K  \frac{n^k_{1}!}{(n^k_{1}-l^k_{1})! l^k_1!}    [1-F_k(t_1)]^{l^k_1}[F_k(t_1)]^{n^k_1-l^k_1}  \frac{n^k_{2}!}{(n^k_{2}-l^k_{2})! l^k_2!} [1-F_k(t_2)]^{l^k_2}[F_k(t_2)]^{n^k_2-l^k_2}\\
&\times \frac{n^k_{12}!}{(n^k_{12}-l^k_{[1]2})!  (l^k_{[1]2}-l^k_{1[2]})! l^k_{1[2]}!}[F_k(t_1)]^{n^k_{12}-l^k_{[1]2}}[F_k(t_2)-F_k(t_1)]^{l^k_{[1]2}-l^k_{1[2]}}[1-F_k(t_2)]^{l^k_{1[2]}}
\end{align}
This probability is similarly derived by $t_1>t_2$ and for $t_1=t_2$.

\section{Joint survival signature of three coherent systems with shared components}
\label{sec.three}

There will be situations where more than two systems may share some components, e.g.\ in networks of systems where a central
server or electricity supply may serve several or even all systems. Further examples can be encountered when systems consisting
of combined hardware and software are considered, with software often shared between different systems. The concept of
joint survival signature, presented in this paper, can be generalized to any number of systems with any kind of component
sharing. To illustrate such a generalization, we briefly consider the case of three systems with a single type of
components, where the systems share some components. Generalization to multiple types of components can be achieved along
the lines presented above.

Consider three coherent systems, $S_1, S_2$ and $S_3$. Let the number of shared components $S_i$ and $S_j$ be denoted by $n_{ij}$,
and let the number of components shared by all three systems be $n_{123}$. Let $n_i$ be the number of components in $S_i$ that are
not shared with any other system. The numbers of components in the systems are
\begin{align*}
&n_1^{*}=n_1+n_{12}+n_{13}+n_{123}\\
&n_2^{*}=n_2+n_{12}+n_{23}+n_{123}\\
&n_3^{*}=n_3+n_{13}+n_{23}+n_{123}
\end{align*}
and the total number of components is $n=n_1+n_2+n_3+n_{12}+n_{13}+n_{23}+n_{123}$.

Let $T_1$, $T_2$ and $T_3$ be the failure times of systems $S_1, S_2$ and $S_3$, respectively, based on 
components with {\it iid} failure times $X_1,X_2, \ldots, X_{n}$ having a common continuous distribution 
function with CDF $F(t)$. Let $l_{[1]2}$ and $l_{1[2]}$ be the number of components out of $n_{12}$ that  
function when, respectively, $S_1$  and $S_2$ is considered. Similarly let $l_{[1]3}$ ($l_{1[3]}$) be the 
number of components out of $n_{13}$ that  function when $S_1$ ($S_3$) is considered, let $l_{[2]3}$ 
($l_{2[3]}$) be the number of components out of $n_{23}$ that  function when $S_2$ ($S_3$) is considered, 
and finally let $l_{[1]23}$,   $l_{1[2]3}$ and $l_{12[3]}$ be the number of components out of $n_{123}$ that  
function when $S_1$, $S_2$ and $S_3$  are  considered, respectively. For ease of notation let 
$\underline{l}=(l_1,l_2,l_3,l_{[1]2},l_{1[2]},l_{[1]3},l_{1[3]},l_{[2]3},l_{2[3]},l_{[1]23},l_{1[2]3},l_{12[3]})$, 
and we denote the summation over all their possible values as $\sum_{\underline{l}}$. The joint survival signature 
for these three systems, $\Phi(\underline{l})$, can be defined as
\begin{equation} 
\Phi(\underline{l})=P(SF_1, SF_2, SF_3 |  \underline{l})
\end{equation} 
where $SF_i$ represents the event that system $S_i$ functions at the time point it is considered.

Generalizing the approach for two systems, as presented above, the joint survival function for these three systems
is derived by
\begin{equation} 
P(T_{1}>t_1,T_{2}>t_2,T_{3}>t_3)=\sum_{\underline{l}}
\Phi(\underline{l})P(C_{\underline{t},\underline{l}})
\end{equation} 
where $P(C_{\underline{t},\underline{l}})$ denotes the probability of the event that the vector $\underline{l}$ describes
precisely the numbers of components, shared or not shared among the systems, that function at the times $\underline{t}=
(t_1, t_2, t_3)$, where $t_i$ is the time at which $S_i$ is considered. For example, for $t_1<t_2<t_3$, which logically
implies $l_{[1]2}\geq l_{1[2]}$, $l_{[1]3}\geq l_{1[3]}$, $l_{[2]3}\geq l_{2[3]}$ and $l_{[1]23}\geq l_{1[2]3}\geq l_{12[3]}$, 
this probability is
{\small
\begin{align}
\nonumber P(C_{\underline{t},\underline{l}})&=P(C^1_{t_1}=l_1,C^2_{t_2}=l_2,C^3_{t_3}=l_3, C^{[1]2}_{t_1}=l_{[1]2},C^{1[2]}_{t_2}=l_{1[2]}, 
C^{[1]3}_{t_1}=l_{[1]3},C^{1[3]}_{t_3}=l_{1[3]},\\
\nonumber &\hspace{1cm}C^{[2]3}_{t_2}=l_{[2]3},C^{2[3]}_{t_3}=l_{2[3]} ,C^{[1]23}_{t_1}=l_{[1]23} ,C^{1[2]3}_{t_2}=l_{1[2]3},C^{12[3]}_{t_3}=l_{12[3]})\\
\nonumber &=\frac{n_{1}!}{(n_{1}-l_{1})! l_1!} [1-F(t_1)]^{l_1}[F(t_1)]^{n_1-l_1}\\
\nonumber &\times \frac{n_{2}!}{(n_{2}-l_{2})! l_1!} [1-F(t_2)]^{l_2}[F(t_2)]^{n_2-l_2}\\
\nonumber &\times \frac{n_{3}!}{(n_{3}-l_{3})! l_1!} [1-F(t_3)]^{l_3}[F(t_3)]^{n_3-l_3}\\
\nonumber &\times \frac{n_{12}!}{(n_{12}-l_{[1]2})! (l_{[1]2}-l_{1[2]})! l_{1[2]}!} [F(t_1)]^{n_{12}-l_{[1]2}}[F(t_2)-F(t_1)]^{l_{[1]2}-l_{1[2]}}[1-F(t_2)]^{l_{1[2]}}\\
\nonumber &\times \frac{n_{13}!}{(n_{13}-l_{[1]3})!  (l_{[1]3}-l_{1[3]})! l_{1[3]}!}[F(t_1)]^{n_{13}-l_{[1]3}}[F(t_3)-F(t_1)]^{l_{[1]3}-l_{1[3]}}[1-F(t_3)]^{l_{1[3]}}\\
\nonumber &\times \frac{n_{23}!}{(n_{23}-l_{[2]3})! (l_{[2]3}-l_{2[3]})!  l_{1[3]}!}[F(t_2)]^{n_{23}-l_{[2]3}}[F(t_3)-F(t_2)]^{l_{[2]3}-l_{2[3]}}[1-F(t_3)]^{l_{1[3]}}\\
\nonumber &\times \frac{n_{123}!}{(n_{123}-l_{[1]23})! (l_{[1]23}-l_{1[2]3})!   (l_{1[2]3}-l_{12[3]})! l_{12[3]}! }\\
&\times [F(t_1)]^{n_{123}-l_{[1]23}}[F(t_2)-F(t_1)]^{l_{[1]23}-l_{1[2]3}}[F(t_3)-F(t_2)]^{l_{1[2]3}-l_{12[3]}}[1-F(t_3)]^{l_{12[3]}}
\end{align}}

Defining this probability similarly for all orderings of $(t_1,t_2,t_3)$ becomes cumbersome with regard to notation, but
the idea will be clear. Further development of this methodology is best done in direct relation to a real-world application,
in order to take specific aspects of the joint systems into account. Extension to multiple types of components, and to more
than three systems, is conceptually trivial following the methodology presented in this paper, but cumbersome with regard to notation.
Suitable computational algorithms for exact calculation or approximations also need to be developed, this is left as an important
topic for future research. The next example briefly illustrates the computational effort required when applying this method in a
naive way to three simple systems with joint components, this particularly serves to emphasize the need for development of 
suitable computational theory and algorithms.

\begin{example}\label{Ex4}{\rm
In addition to the two systems in Example \ref{Ex2}, we consider a third system as in Figure \ref{fig.ex4}.  
To illustrate the complexity of the approach in this paper for more than two systems, we show the effort required to compute the survival signature restricted to considering all systems at the same moment of time. Both systems 1 and 2 have 3 components while system 3 has 4 components, and we have $n_1=n_2=0$, $n_3=1$ (E), $n_{12}=1$  (B), $n_{13}=1$  (C), $n_{23}=1$  (D), $n_{123}=1$  (A).  As in this example we have at most one shared component among the systems, thus $l_1=l_2=0$, $l_3\in\{0,1\}$, $l^*_{12}=l_{[1]2}=l_{1[2]}\in\{0,1\}$, $l^*_{13}=l_{[1]3}=l_{1[3]}\in\{0,1\}$, $l^*_{23}=l_{[2]3}=l_{2[3]}\in\{0,1\}$, $l^*_{123}=l_{[1]23}=l_{1[2]3}=l_{12[3]}\in\{0,1\}$. As $l_1=l_2=0$ there are $2^{10}$ possibilities we need to consider, the first column in Table \ref{tab.ex4} shows the corresponding number of these possibilities out of $2^{10}$. For example, the first three rows suggest that at least one of these  systems fails when either component B or D fails, and this  counts for $3\times 2^8=768$ out of $2^{10}=1024$ possibilities. Clearly, this already requires many combinations to be considered, and the combinatorics increase enormously when also considering the systems at different moments in time. 

  \input{fig.ex4} 
 
 \begin{table}[H]
  \centering\small
    \begin{tabular}{ccccccc}
&  \multicolumn{1}{c}{B} & \multicolumn{1}{c}{D} & \multicolumn{1}{c}{A} & \multicolumn{1}{c}{C} & \multicolumn{1}{c}{E}& \\
$\#$ &   \multicolumn{1}{c}{$l^*_{12}\in \{0,1\}$} & \multicolumn{1}{c}{$l^*_{23}\in \{0,1\}$} & \multicolumn{1}{c}{$l^*_{123}\in \{0,1\}$}  &  \multicolumn{1}{c}{$l^*_{13}\in \{0,1\}$} &\multicolumn{1}{c}{$l_{3}\in \{0,1\}$}& \multicolumn{1}{c}{All systems function} \\[1ex]
\hline
$2^8$	&	0	&	0	&	0/1	&	0/1	&	0/1	&	No\\
$2^8$	&	0	&	1	&	0/1	&	0/1	&	0/1	&	No\\
$2^8$	&	1	&	0	&	0/1	&	0/1	&	0/1	&	No\\
$2^7$	&	1	&	1	&	1	&	0/1	&	0/1	&	Yes\\
$2^6$	&	1	&	1	&	0	&	0	&	0/1	&	No\\
$2^5$	&	1	&	1	&	0	&	1	&	1	&	Yes\\
$2^5$	&	1	&	1	&	0	&	1	&	0	&	No\\
\hline
    \end{tabular}
  \label{tab.ex4}
\caption{Three systems with one type of components, Example \ref{Ex4}}
\end{table}

}\end{example}

\section{Concluding Remarks} 
\label{concl}

In this paper we have introduced the concept of joint survival signature for two systems with shared components. 
First we considered the case when we have only one type of component, then we extended that for multiple types of 
components. We showed how this can be generalised for more than two systems with one or multiple types of shared components.  
We have also presented how the joint survival signature can be used to derive marginal and conditional survival functions.

It is possible to derive variations to the presented joint survival signature for the case with multiple systems sharing
components, and for some specific scenarios reduced versions of the joint survival signature presented here may be
sufficient. For example, one could only take into account the total numbers of functioning components of each type
per system, and use the theorem of total probability and assumed exchangeability of the failure times of components of
each type to relate this to our survival signature and to inferences on the systems' reliability. These suggest interesting
directions for future research, which will be particularly useful if motivated by practical applications.

A crucial consideration is how the joint survival signature can be computed. This is mostly left as an important topic for
future research, in the examples in this paper only small systems are considered for which all combinations of functioning
components and the corresponding state of the two systems are easily checked. One can use the classical theory of minimal cut and
path sets for small systems, but it will be important to develop algorithms to compute the joint survival signature for more
complex systems, in particular where there are multiple types of components and possibly more than two systems being
considered. A main advantage of the joint survival signature for coherent systems is that it is a non-decreasing function
of all its inputs, hence it is straightforward to derive bounds if one only computes the function for a limited number of
inputs. The question which inputs to focus on in order to derive useful bounds with relatively little computation time 
is also interesting for future research. Throughout, a major advantage of the joint survival signature over the full 
structure function is that it requires substantially less storage, which particularly for large systems with relatively
few components types can be very important. While computing survival signatures may be cumbersome, and may require approximations
to be developed, for example by simulations of the system for certain inputs, the computation is only required once for a 
system with a fixed structure.

In Section \ref{sec.sur.sign}, a brief discussion was provided of recent developments based on the single system
survival signature for research and application. All these topics are also of great interest based on the joint survival
signature, where for example consideration of component importance brings novel aspects to the literature as it is likely
that components shared between multiple systems are more important when all systems are considered than they may be for
a single system.

\section*{Acknowledgements}

We thank three reviewers for their support of the paper and useful comments that have led to improved presentation.
The results presented in this paper were mostly achieved during Prof.\ Balakrishnan's visit to Durham University in 
November 2018, funded by a Durham University Global Engagement Travel Grant, which is gratefully acknowledged. 
The work has been presented by the first author at the UK Reliability Meeting, Durham, 1-3 April 2019 
(\url{http://www.maths.dur.ac.uk/stats/uk-reliability}).



\end{document}

%% file: fig.ex1a.tex
\begin{figure}[t]
\centering\setlength{\unitlength}{1mm}
\begin{picture}(100,30)
  \setlength\fboxsep{0pt}  
\put(5,35){System 1}

\put(3,5){\colorbox{gray!20}{\framebox(7,7)[cc]{B}}}
\put(3,15){\colorbox{gray!20}{\framebox(7,7)[cc]{A}}}
\put(3,25){\framebox(7,7)[cc]{D}}

\put(0,9){\line(0,1){20}}
\put(13,9){\line(0,1){20}}

\put(0,9){\line(1,0){3}}
\put(0,29){\line(1,0){3}}
\put(-5,19){\line(1,0){5}}

\put(10,9){\line(1,0){3}}
\put(10,29){\line(1,0){3}}
\put(13,19){\line(1,0){7}}

\put(0,19){\line(1,0){3}}
\put(10,19){\line(1,0){3}}


\put(23,11){\colorbox{gray!20}{\framebox(7,7)[cc]{C}}}

\put(23,21){\framebox(7,7)[cc]{E}}

\put(20,14){\line(0,1){11}}
\put(33,14){\line(0,1){11}}

\put(20,14){\line(1,0){3}}
\put(20,25){\line(1,0){3}}

\put(30,14){\line(1,0){3}}
\put(30,25){\line(1,0){3}}

\put(33,19){\line(1,0){5}}


\put(65,35){System 2}

\put(83,5){\colorbox{gray!20}{\framebox(7,7)[cc]{C}}}
\put(83,15){\colorbox{gray!20}{\framebox(7,7)[cc]{B}}}
\put(83,25){\framebox(7,7)[cc]{G}}

\put(80,9){\line(0,1){20}}
\put(93,9){\line(0,1){20}}

\put(80,9){\line(1,0){3}}
\put(80,29){\line(1,0){3}}
\put(55,19){\line(1,0){5}} 
\put(75,19){\line(1,0){5}}

\put(90,9){\line(1,0){3}}
\put(90,29){\line(1,0){3}}
\put(93,19){\line(1,0){7}}

\put(80,19){\line(1,0){3}}
\put(90,19){\line(1,0){3}}


\put(63,11){\colorbox{gray!20}{\framebox(7,7)[cc]{A}}}

\put(63,21){\framebox(7,7)[cc]{F}}

\put(60,14){\line(0,1){11}}
\put(73,14){\line(0,1){11}}

\put(60,14){\line(1,0){3}}
\put(60,25){\line(1,0){3}}

\put(70,14){\line(1,0){3}}
\put(70,25){\line(1,0){3}}

\put(73,19){\line(1,0){5}}

\end{picture}
\caption{Two systems with one type of components, Example \ref{fig.ex1a}}
\label{fig.ex1a}
\end{figure}

%% file: fig.ex2.tex
\begin{figure}[t]

\centering\setlength{\unitlength}{1mm}
\begin{picture}(100,30)
  \setlength\fboxsep{0pt}  
\put(0,28){System 1}

\put(3,15){\colorbox{gray!20}{\framebox(7,7)[cc]{A}}}
\put(3,5){\colorbox{gray!20}{\framebox(7,7)[cc]{B}}}
\put(20,15){\framebox(7,7)[cc]{C}}

\put(-5,14){\line(1,0){5}}
\put(30,14){\line(1,0){5}}

\put(0,9){\line(0,1){10}}
\put(30,9){\line(0,1){10}}

\put(0,9){\line(1,0){3}}

\put(10,9){\line(1,0){20}}

\put(0,19){\line(1,0){3}}
\put(10,19){\line(1,0){10}}
\put(27,19){\line(1,0){3}}


\put(50,28){System 2}

\put(73,9){\colorbox{gray!20}{\framebox(7,7)[cc]{B}}}

\put(45,13){\line(1,0){5}} 
\put(65,13){\line(1,0){5}}

\put(83,13){\line(1,0){3}}

\put(70,13){\line(1,0){3}}
\put(80,13){\line(1,0){3}}

\put(53,15){\colorbox{gray!20}{\framebox(7,7)[cc]{A}}}

\put(53,5){\framebox(7,7)[cc]{D}}

\put(50,8){\line(0,1){11}}
\put(63,8){\line(0,1){11}}

\put(50,8){\line(1,0){3}}
\put(50,19){\line(1,0){3}}

\put(60,8){\line(1,0){3}}
\put(60,19){\line(1,0){3}}

\put(63,13){\line(1,0){5}}

\end{picture}
\label{fig.ex2}
\caption{Two systems with one type of components, Example \ref{Ex2}}
\end{figure}

%% file: fig.ex4.tex
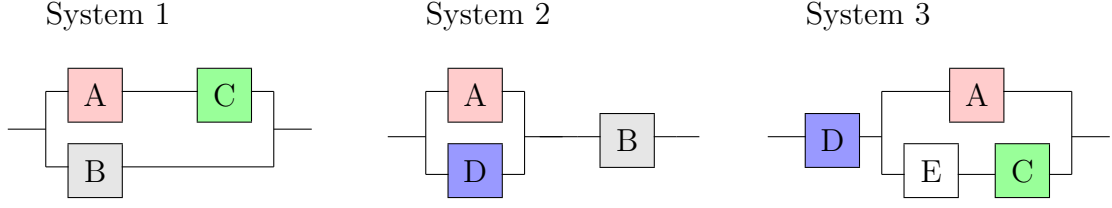
\begin{figure}[t]

\centering\setlength{\unitlength}{1mm}
\begin{picture}(135,30)
  \setlength\fboxsep{0pt}  
\put(0,28){System 1}

\put(3,15){\colorbox{red!20}{\framebox(7,7)[cc]{A}}}
\put(3,5){\colorbox{gray!20}{\framebox(7,7)[cc]{B}}}
\put(20,15){\colorbox{green!40}{\framebox(7,7)[cc]{C}}}

\put(-5,14){\line(1,0){5}}
\put(30,14){\line(1,0){5}}

\put(0,9){\line(0,1){10}}
\put(30,9){\line(0,1){10}}

\put(0,9){\line(1,0){3}}

\put(10,9){\line(1,0){20}}

\put(0,19){\line(1,0){3}}
\put(10,19){\line(1,0){10}}
\put(27,19){\line(1,0){3}}


\put(50,28){System 2}

\put(73,9){\colorbox{gray!20}{\framebox(7,7)[cc]{B}}}

\put(45,13){\line(1,0){5}} 
\put(65,13){\line(1,0){5}}

\put(83,13){\line(1,0){3}}

\put(70,13){\line(1,0){3}}
\put(80,13){\line(1,0){3}}

\put(53,15){\colorbox{red!20}{\framebox(7,7)[cc]{A}}}

\put(53,5){\colorbox{blue!40}{\framebox(7,7)[cc]{D}}}

\put(50,8){\line(0,1){11}}
\put(63,8){\line(0,1){11}}

\put(50,8){\line(1,0){3}}
\put(50,19){\line(1,0){3}}

\put(60,8){\line(1,0){3}}
\put(60,19){\line(1,0){3}}

\put(63,13){\line(1,0){5}}


\put(100,28){System 3}

\put(100,9){\colorbox{blue!40}{\framebox(7,7)[cc]{D}}}

\put(95,13){\line(1,0){5}} 
\put(107,13){\line(1,0){3}}

\put(120,8){\line(1,0){5}}


\put(119,15){\colorbox{red!20}{\framebox(7,7)[cc]{A}}}

\put(113,5){\framebox(7,7)[cc]{E}}
\put(125,5){\colorbox{green!40}{\framebox(7,7)[cc]{C}}}

\put(110,8){\line(0,1){11}}
\put(135,8){\line(0,1){11}}

\put(132,8){\line(1,0){3}}
\put(126,19){\line(1,0){9}}

\put(110,8){\line(1,0){3}}
\put(110,19){\line(1,0){9}}

\put(135,13){\line(1,0){5}}

\end{picture}
\caption{Three systems with one type of components, Example \ref{Ex4}}
\label{fig.ex4}
\end{figure}